\documentstyle[12pt]{article}

\newfont{\frak}{eufm10 scaled\magstep1}
\newfont{\extra}{msbm10 scaled\magstep1}

\newcommand{\fra}[1]{\mbox{\frak #1}}
\newcommand{\extr}[1]{\mbox{\extra #1}}

\newcommand{\sect}[1]{\setcounter{equation}{0}\section{#1}}
\newcommand{\subsect}[1]{\subsection{#1}}

\renewcommand{\theequation}
  {\arabic{section}.\arabic{equation}}

\def\ie{{\it i.e.}}
\def\eg{{\it e.g.}}

\def\beq{\begin{equation}}
\def\eeq{\end{equation}}
\def\ba{\begin{array}}
\def\ea{\end{array}}

\def\U{{\cal U}}
\def\Fun{\hbox{Fun}}
\def\co{\Delta}                      

\def\a{\alpha}
\def\b{{\beta}}

\def\bbeta{\bar{\beta}}

\def\v{\varepsilon}

\def\k{\omega}                

\def\J{{\extr J}}            
\def\P{{\extr P}}            

\def\R{{\extr R}}  
\def\T{{\extr T}}  
\def\X{{\extr X}}  
\def\W{{\extr W}}  

\def\C{\mbox{\ C}}      
\def\S{\mbox{\ S}} 

\def\s{{{so}}}           
\def\is{{{iso}}}           

\def\kin{\k_1,\k_2,\dots,\k_N}
\def\koiin{0,\k_2,\dots,\k_N}
\def\kiin{\k_2,\dots,\k_N}

\def\pardeform{\lambda}

\def\Gq{\Gamma^{(m)}_\pardeform}
\def\Gc{\Gamma^{(m)}}

\catcode`\@=11
\font\tenmsa=msam10
\font\sevenmsa=msam7
\font\fivemsa=msam5
\font\tenmsb=msbm10
\font\sevenmsb=msbm7
\font\fivemsb=msbm5
\newfam\msafam
\newfam\msbfam
\textfont\msafam=\tenmsa  \scriptfont\msafam=\sevenmsa
  \scriptscriptfont\msafam=\fivemsa
\textfont\msbfam=\tenmsb  \scriptfont\msbfam=\sevenmsb
  \scriptscriptfont\msbfam=\fivemsb

\def\hexnumber@#1{\ifnum#1<10 \number#1\else
 \ifnum#1=10 A\else\ifnum#1=11 B\else\ifnum#1=12 C\else
 \ifnum#1=13 D\else\ifnum#1=14 E\else\ifnum#1=15 F\fi\fi\fi\fi\fi\fi\fi}

\def\msa@{\hexnumber@\msafam}
\def\msb@{\hexnumber@\msbfam}
\mathchardef\blacktriangleright="3\msa@49
\mathchardef\blacktriangleleft="3\msa@4A
\catcode`\@=12

\def\bicross{\triangleright\!\!\!\blacktriangleleft}
\def\rimo{\triangleright\!\!\!<}
\def\leco{>\!\!\blacktriangleleft}

\def\LR{\blacktriangleright\!\!\!\triangleleft}

\def\RRR{\bar\triangleright}

\def\AA{H}
\def\HH{A}
\def\KK{K}
\def\aa{h}
\def\bb{g}
\def\hh{a}
\def\gg{b}

\begin{document}

\rightline{DAMTP 96--100}
\rightline{To appear in J. Phys. {\bf A}}
\vspace{1.5cm}

\begin{center} 
{\large{\bf{GRADED CONTRACTIONS AND BICROSSPRODUCT STRUCTURE}}}
{\large{\bf{OF DEFORMED INHOMOGENEOUS ALGEBRAS 
}}}
\end{center}  

\bigskip\bigskip

\begin{center} 
J. A. de Azc\'arraga$^{1}$ 
\footnote{St. John's College Overseas Visiting Scholar.}
\footnote{On sabbatical (J.A.) leave and on leave of absence (J.C.P.B.)
from Departamento de F\'{\i}sica Te\'orica and IFIC, Centro Mixto Univ. 
de Valencia-CSIC, E--46100 Burjassot (Valencia), Spain. E-mails: 
j.azcarraga@damtp.cam.ac.uk (azcarrag@evalvx.ific.uv.es), 
pbueno@lie.ific.uv.es.},
M.A. del Olmo$^2$  \footnotemark[4],  
J.C. P\'erez Bueno$^{1}$ \footnotemark[3], 
and M. Santander$^2$ 
\footnote{e-mails: olmo@cpd.uva.es, santander@cpd.uva.es.}
\end{center}

\begin{center} {\it $^1$ Department of Applied Mathematics and Theoretical 
Physics, \\
Silver St., Cambridge CB3 9EW, UK}
\end{center}

\begin{center}{\it{$^2$ Departamento de F\'{\i}sica Te\'orica,
Universidad de Valladolid}\\ E--47011, Valladolid, Spain} 
\end{center}


\begin{abstract}
A family of deformed Hopf algebras corresponding to the classical
maximal isometry algebras of zero-curvature $N$-dimensional spaces (the
inhomogeneous algebras ${iso}(p,q), \ p+q=N,$  as well as some of their
contractions) are shown to have a bicrossproduct structure. This is done for
both the algebra and, in a low-dimensional example, for the (dual) group
aspects of the deformation.
\end{abstract}

\vfill\eject
\sect{Introduction}

The procedure to deform simple algebras and groups
was established by Drinfel'd
\cite{Dri}, Jimbo \cite{Ji} and Faddeev, Reshetikhin and Takhtajan \cite{FRT}. 
The algorithm, which leads to the
so-called  `quantum' algebras, does not cover, however, the case of
non-semisimple algebras. Since the contraction process leads to
inhomogeneous algebras by starting from  simple ones, it is natural
to use it as a way to deform inhomogeneous Lie  (\ie, `classical'
or undeformed) algebras. This path of extending the classical idea
of the Lie algebra contraction to the case of deformed  algebras
was proposed by Celeghini {\em et al.} \cite{CGST}.  The basic
requirement to define a deformed inhomogeneous algebra  is
the commutativity of the processes  of contraction and
deformation: when considering a simple algebra and one of their
inhomogeneous contractions, both at classical and deformed
levels, the deformation of the contracted inhomogeneous Lie algebras 
should coincide with the contraction of the deformed simple algebra.
This commutativity is not always guaranteed, and in general requires 
\cite{CGST} a redefinition of the deformation parameter $q$ in 
terms of the contraction parameter and the  new deformation one, so that
$q$  is not a passive element in the contraction. 
This was used, for instance,  to obtain the $\kappa$-Poincar\'e algebra
\cite{LNRT}, for which the  deformation parameter
$\kappa$ has dimensions of inverse of length.  

The concept of contraction of Lie
algebras (or groups) was discussed in the early fifties by \.In\"on\"u 
and Wigner \cite{IW} (see also \cite{Saletan}).   
The idea of group contraction itself arose in the group analysis
of the non-relativistic limit, and its applications to
mathematical physics problems have been very fruitful. 
The study of details
behind this procedure unveils interesting mathematical structures,
which in many important cases are linked to physical properties. 
In particular, the contraction process may increase the group cohomology
\cite{AA} (see also \cite{AHPS}), 
as it is the case in the standard non-relativistic limit. 
Several attempts have been made to systematise the study of
contractions 
and recently a new approach has
been put forward by Moody, Montigny and Patera \cite{dMPMP},
under the name of graded contractions. 
The key idea there is to preserve
a given grading of the original Lie algebra. This condition may fit
neatly with physical requirements and is automatically
satisfied in the simplest case of the \.In\"on\"u--Wigner contractions,
which correspond to the simplest ${\extr Z}_2$-grading.

A class of Lie algebras describing a whole family of
contractions is the so-called orthogonal Cayley-Klein (CK) algebras.  The
name is due to historical reasons:  these are the Lie algebras of the
motion groups of real spaces with a 
projective metric \cite{SommerYRY} (see also \cite{SHO}). 
The same family appears as a natural subset of all ${\extr Z}_2^{\otimes
N}$-graded contractions which can be obtained from $\s(N+1)$ \cite{HMOS}. 
And furthermore, among orthogonal CK
algebras we find not only all simple pseudo-orthogonal algebras,
but many non-semisimple algebras of physical importance, as the
kinematical Poincar\'e and Galilei algebras in $(N-1,1)$ dimensions, the
Euclidean algebra in $N$ dimensions, etc. The CK scheme
does not deal with a single Lie algebra, but with a whole
family of them simultaneously,
each of which is parametrised by a set of real
numbers with a well-defined geometrical and physical significance. 
The main point to be stressed is the ability of this kind of approach to
describe some properties of many Lie algebras in a single unified form.
This is possible as the Lie algebras in the CK family, though not
simple, are `very near' to the simple ones, and many structural
properties of the simple algebras, when suitably reformulated, still
survive for the CK algebras.

It is possible to give deformations of algebras in the CK
family; naturally enough these will be 
said to belong to the CK family of Hopf `quantum'  algebras. 
In \cite{BHOSab} deformations of the
enveloping algebras of all algebras in the CK family of 
$\s(p,q),\ p+q=3, 4$ were given. For higher dimensions, \ie\ for 
algebras in the family of $\s(p,q),\ (p+q=N+1)$ with $N>3$, a
quantum deformation of the general parent member of the CK family 
is still not known, yet there exists a scheme of quantum deformations
encompassing all  motion algebras of flat affine spaces in $N$
dimensions, which include the ordinary inhomogeneous 
$\is(p,q),\ (p+q=N)$ 
\cite{BHOScd}. This scheme provides a Hopf algebra
deformation for each algebra in the family. Some of its members 
are physically relevant non-semisimple algebras, and include as
particular cases most of the deformations of these algebras found in the
literature.   

An important fact in quantum algebra/group theory is the (co)existence
of two closely linked algebraic structures: the algebra (as expressed by the 
commutators or the commuting properties of the algebra 
of functions on the group) and the coalgebra (as given by the coproduct). 
Most of the complications found when doing quantum contractions can
be traced to the need of dealing simultaneously with these two aspects.
For instance, a naive contraction might lead to divergences either in
the coproduct or in the $R$--matrix \cite{CGST,BGHOS}.
One of the main motivations behind the CK scheme was to
be able to describe at the same time a family of algebras, including
some simple and some contracted algebras, in such a way that the
possible origin of divergences under contractions is clearly seen
and controlled. 

In this paper we address a specific problem where the advantages of a
CK type scheme are  exhibited.
In the classical case, an \.In\"on\"u--Wigner (IW) contraction
of a simple algebra leads to a non-semisimple one which is the
semidirect sum of an abelian algebra and the preserved subalgebra of
the original algebra with respect which the contraction was made.
{\it All} IW  contractions of simple algebras have a semidirect
structure. 
It is then natural to ask: is there a similar pattern for the contracted 
deformations \ie,
for the Hopf algebra deformations of contracted simple Lie algebras? 
The analogue of the semidirect product is an example of the bicrossproduct 
of Hopf algebras, introduced by Majid \cite{MajMajid} 
(see also \cite{SingMoln,bicII}).
The aim of this paper is to show that all deformed algebras in
the affine\footnote{We use the word `affine' in the sense of inhomogeneous. 
Not all deformed inhomogeneous groups have a bicrossproduct structure; this 
is, for instance, the case of ${\cal U}_q({\cal E}(2))$ 
as discussed in \cite{APb}.} 
CK family $\is_{\kiin}(N)$ have indeed a
bicrossproduct structure, as is the case of the $\kappa$-Poincar\'e 
\cite{KPoinBicros}.  
This result opens the possibility of
recovering more easily the deformed dual groups
$\Fun_q(ISO_{\kiin}(N))$ by using the dual bicrossproduct
`group-like' expressions (see \cite{APb} for some group-like (rather
than algebra-like) examples  of this construction).
Classically, the $\is_{\kiin}(N)$ family includes all inhomogeneous
Lie algebras $\is (p,q),\ (p+q=N)$, so we will refer loosely to the aim of 
the paper as showing the bicrossproduct structure of deformed inhomogeneous
groups. It should be kept in mind, however, that we are referring to a
specific deformation, and that there exist examples 
(see \cite{APb}) where 
a contraction of a deformed algebra has no bicrossproduct structure. 

The paper is organised as follows. In Section II we  
briefly describe the classical Cayley--Klein algebras and  
present a discussion on  contractions and 
dimensional analysis since this is relevant for the assignment of physical 
dimensions to the deformation parameters. 
In Sec. III we give the explicit expressions for their
$q$--deformations. The bicrossproduct structure of these $q$--deformed
Cayley--Klein Hopf algebras is shown in Section IV.  Examples of this structure
for physically interesting algebras are presented in Section V. 
In Section VI we show, as an example,  how to obtain the (dual) 
group deformation in the case of lowest dimension $N=2$. 
Some conclusions close the paper. 


\sect{ Affine Cayley--Klein Lie algebras and dimensional analysis}

\subsect{The CK scheme of geometries and Lie algebras}

\bigskip
The complete family of the $so(N+1)$ CK algebras is a set of real Lie algebras 
of dimension $(N+1)N/2$, characterised by
$N$ real parameters $(\k_1,\k_2,\dots,\k_N)$ \cite{SHO}.  
This family appears, \eg, as a
natural subfamily \cite{HMOS} of all the graded
contractions from the Lie algebra $so(N+1)$ \cite{HS} corresponding to a
${\extr Z}_2^{\otimes N}$ grading of $so(N+1)$, 
and its elements will be denoted
$\s_{\kin} (N+1)$; 
in particular, 
$\s_{1, 1, \dots, 1} (N+1) \equiv \s(N+1)$. 
In terms of a basis of
$\s_{\kin} (N+1)$ adapted to the grading, $\{ \J_{ab};\ a<b,\ a,b= 0,
1, \dots, N\}$,  
this family of algebras is defined by
\beq
[\J_{ab},\J_{ac}] =  \k_{ab} \J_{bc}, \quad
[\J_{ab},\J_{bc}] = -\J_{ac} ,  \quad
[\J_{ac},\J_{bc}] = \k_{bc}\J_{ab},\label{aa}
\eeq
where now $a<b<c, \ a,b,c =0,1,\dots ,N$,
$\k_{ab}:=\k_{a+1}\k_{a+2}\ldots\k_b=\prod_{l=a+1}^{b}\k_{l}$
(thus, $\k_{ab}\k_{bc}=\k_{ac}$) and $[\J_{ab}, \J_{cd}]=0$ if the four indices 
are different. 
By simple rescaling of the generators, all the numerical values of the
constants $\k_i$ may be brought to one of the  values $1, 0, -1$, hence the
complete CK family contains $3^N$ algebras which are  different as graded
contractions, even if 
some of them may still be isomorphic.

When all the $\k_i$ are non-zero but some of them
are negative, the algebra $\s_{\kin} (N+1)$ is isomorphic to a certain
pseudo-orthogonal algebra $\s (p,q)$ $(p+q=N+1,\, p\ge q > 0)$. 
If all the $\k_i$ are non-zero  we can introduce also
$\J_{ba}, (a<b)$ by $\J_{ba}:= -\frac{1}{\k_{ab}}\J_{ab}$ and $\J_{aa}
:=0$, so that the commutation relations can be written in the familiar form: 
\beq
[\J_{ij},\J_{lm}]=\delta_{im}\J_{lj}-\delta_{jl}\J_{im}+\delta_{jm}\k_{lm} \J_{il}
+\delta_{il}\k_{ij} \J_{jm}.
\label{etiqueta}
\eeq
If, however, some constant(s) $\k_i=0$, the
algebras (\ref{aa}) become inhomogeneous and correspond to algebras that
are obtained from $\s (p,q)$
through a sequence of \.In\"on\"u-Wigner (IW) contractions. 
To describe them let us denote ${\fra h}^{(m)}$ ($m=1, \dots, N$) the 
subalgebra generated by the $\J_{ab}$ ($a<b$) for which
$a,b$  satisfy either $b< m$ {\it or} $a \geq m$. 
A complement for ${\fra h}^{(m)}$ is the vector
subspace ${\fra p}^{(m)}$ (not always a subalgebra) spanned by the
elements $\J_{ab}$ with $a < m$ 
{\it and} $b \geq m$. The decomposition
$\s_{\kin} (N+1) = {\fra p}^{(m)} \oplus {\fra h}^{(m)}$ is in fact a
Cartan-like decomposition, and there exists an involutive automorphism
of the Lie algebra $\s_{\kin} (N+1)$ with  ${\fra p}^{(m)}$ and 
${\fra h}^{(m)}$ as the 
anti-invariant and invariant subspaces. 
The structure of the subalgebra ${\fra h}^{(m)}$ and of the vector subspace 
${\fra p}^{(m)}$ of the Lie
algebra $\s_{\kin} (N+1)$ can be graphically displayed by arranging the 
generators of $so(N+1)$ in the form of a  triangle

\bigskip

\noindent\hskip 1.8truecm 
\beq\begin{tabular}{cccc|cccc}
$J_{01} $&$ J_{02} $&$\ldots$&
$J_{0\, (m-1)} $&
 $J_{0m}$&$J_{0\, (m+1)}$&
$\ldots$&$J_{0N}$\\
 &$ J_{12} $&$\ldots$& $J_{1\, (m-1)} $& $J_{1m}$&$J_{1\, (m+1)}$&
$\ldots$&$J_{1N}$\\
 &&$\ddots $&$\vdots$&  $\vdots$&$\vdots$&
$ $&$\vdots$ \\
 & &$ $&$J_{(m-2)\,(m-1)}$&  $J_{(m-2)\,m}$&$J_{(m-2)\,(m+1)}$&
$\ldots$&$J_{(m-2)\,N}$\\
 & &$ $& &  $J_{(m-1)\,m}$&$J_{(m-1)\,(m+1)}$&
$\ldots$&$J_{(m-1)\,N}$ \label{triangulo}\\
\cline{5-8}
 & &$ $&\multicolumn{1}{c}{\,}&    $ $&$J_{m\,(m+1)}$&
$\ldots$&$J_{m  N}$\\
 & &$ $&\multicolumn{1}{c}{\,}& $ $ &$ $&
$\ddots$&$\vdots$\\
 & &$ $&\multicolumn{1}{c}{\,}& $ $&$ $& 
$ $&$J_{N-1\,N}$\\
\end{tabular}
\eeq

\bigskip
\noindent
We see that the generators which 
span the subspace ${\fra p}^{(m)}$ are the
$m(N+1-m)$ generators in the rectangle 
determined by the corner $J_{(m-1)m}$. 
The triangles at its left and below correspond to the subalgebras
$\s_{\k_1\ldots \k_{m-1}}(m)$ and 
$\s_{\k_{m+1}\ldots \k_N}(N+1-m)$ respectively, 
the direct sum of which is the subalgebra 
${\fra h}^{(m)}$. The  subspace
${\fra p}^{(m)}$ 
corresponding to the $\k_m$-rectangle in the
diagram, can be identified 
with the Lie algebra quotient space 
$so_{\kin} (N+1) / \s_{\k_1\ldots \k_{m-1}}(m) \oplus \s_{\k_{m+1}\ldots
\k_N}(N+1-m)$.

For each decomposition $\s_{\kin} (N+1) = {\fra p}^{(m)} \oplus {\fra h}^{(m)}, 
\ m=1,\dots,N$, 
there is a possible IW contraction,  denoted $\Gc$, to be performed on the
algebra $\s_{\kin} (N+1)$. Specifically, if we denote the generators of 
the standard $\s_{\kin} (N+1)$ algebra by $\extr X$, the IW contraction
$\Gc$ of $\s_{\kin} (N+1)$ is given by the $\epsilon \to 0$
limit of the replacements 
\beq 
\begin{array}{l}
\Gc (\extr X)\equiv \extr X' = \left\{\begin{array}{lc}
\  \extr X \quad &\mbox \ \mbox{if}\ \extr X\in \fra h^{(m)}\\
\epsilon   \extr X \quad &\ \mbox{if}\ \extr X \in \fra p^{(m)}
\end{array}
\right. ,\qquad   m=1,\dots,N, 
\label{ae}
\end{array}\eeq
Under the contraction $\Gc$, the algebra $so_{\kin}(N+1)$ goes to another
algebra in the CK family with  the same values of the $\k_i$ constants except
for $\k_{m} =0$. 
Thus, in the triangular arrangement of
generators, the
$N$ possible IW contractions correspond to the
$N$ different rectangles that can be selected inside the large triangle.
These rectangles are completely abelianised by the contractions, while the
commutators with
one or two generators outside ${\fra p}^{(m)}$ remain unchanged. 
As an example, the contraction given by
(\ref{ae}) with $m=1$ and starting from a $so(p,q)$ algebra, where all $\k_i$ are
different from zero, corresponds to the limit $\epsilon \to 0$ of 
$J_{0i}\mapsto J_{0i}'=\epsilon   J_{0i}$, $J_{ij}=J_{ij}'\ (i\ne 0)\,,\,
J_{0i},J_{ij}\in so(p,q)$. This leads to
$[J_{0i}',J_{0j}']=\pm \epsilon^2 J_{ij}'$ and hence
$J_{0i}'$ 
$i=1,\ldots,N$ determines the abelian $N$-dimensional ideal ${\fra p}^{(1)}$.

Let us now consider the homogeneous space ${\cal S} \equiv SO_{\kin} (N+1) /
SO_{\kiin}(N)$, where $SO_{\kiin}(N)$ is the subgroup generated by
the subalgebra $\fra h^{(m)}$ with $m=1$. This space has an invariant
canonical connection, and a hierarchy of metrics, coming after suitable rescalings
from the Cartan-Killing form  in the algebra
$so_{\kin}(N+1)$.  When 
the constants $\k_{2},\dots, \k_{N}$  are
different from zero, then the `main' metric is non-degenerate, the invariant
canonical connection turns out to be the corresponding Levi-Civita metric
connection, and the space ${\cal S}$ has a curvature which is constant and 
equal to $\k_1$. 

In the particular case $( \kin ) = (\k_1, 1, \dots, 1 )$ the space
${\cal S}$ reduces to the Riemannian space (positive definite metric) of 
constant curvature $\k_1$ and dimension $N$. 
When $\k_1=0$ the algebras $\s_{\koiin} (N+1)$, can be realized as algebras  
of groups of affine transformations on $\R ^N$  \cite{SHO}; in this case we 
shall rename the generators as 
$\{\P_{i}:=
\J_{0i},\
\J_{ij};\ i<j,\ i,j= 1, \dots, N\}$, the new names stressing the role of
$\P_i$ as generating translations and of $\J_{ij}$ as the generators of
rotations around the origin of the space. 
Each subalgebra ${\fra h}^{(m)},\  m=1, \dots ,N$, is spanned by the set of 
generators $\{\P_i,\J_{ij},\ i,j=1,\dots,m-1;\ \J_{kl},  \ k,l=m,\dots,N \}$, 
and hence the collection of subalgebras ${\fra h}^{(m)}$ can be clearly 
identified in the space
${\cal S}$ as the isotropy subalgebras of a point (for $m=1$), of a line (for
$m=2$), \dots , of a hyperplane (for $m=N$).  

The non-zero Lie brackets of $\s_{\koiin} (N+1)$ are given by 
\beq
[\J_{ij},\P_i] = \P_j,  \quad
[\J_{ij},\P_j] = \!-\!\k_{ij}\,\P_i,  	\quad  
[\J_{ij},\J_{ik}] =  \k_{ij} \J_{jk}, \quad
[\J_{ij},\J_{jk}] = -\J_{ik} ,  \quad
[\J_{ik},\J_{jk}] = \k_{jk}\J_{ij},
\label{ab}
\eeq
where the indices $i, j, k =1,\ldots, N$ are always assumed to be ordered,
$i<j<k$.  
Note in particular that all translation generators commute
(as witnessing the zero curvature). 
It will be convenient to denote
this $\s_{\koiin}(N+1)$ Lie algebra by $\is_{\kiin}(N)$ , and by 
$ISO_{\kiin}(N)$ the
corresponding group. There are $3^{N-1}$  different
$N$--dimensional affine CK geometries, and from relations (\ref{ab})
it is clear that the groups $ISO_{\kiin} (N)$ have a semidirect
product structure 
\beq
ISO_{\kiin} (N) = SO_{\kiin} (N)\odot T_N\quad,\label{ac}
\eeq
where $T_N$ is the abelian subgroup generated by $\{ \P_i\ ; \
i=1,\dots,N
\}$ (in the case $\k_1=0$, this abelian subgroup can be identified with
the CK homogeneous space $\cal S$ itself) and
$SO_{\kiin} (N)$ is a general CK group with $N-1$ constants $\k_i$,  
generated by $\{ \J_{ij}\ ; \ 
i,j=1,\dots,N\} $. The `main' 
metric which is kept invariant by
the action of this group is described by the quadratic form given
by a matrix with diagonal entries 
$(1, \ \k_2, \ \k_2 \k_3,\ \dots, \ \k_2 \cdots \k_N)$. Among these
inhomogeneous groups we can recognise the Euclidean group in $N$
dimensions for which 
$(\kin)=(0,1, 1, \dots, 1 )$, the Poincar\'e group in $(N-1,1)$ 
dimensions (appearing several times in the CK affine scheme as
\eg, for $(0,-(1/c^2), 1,\dots, 1 )$) 
or the Galilei group in $(N-1,1)$ dimensions
which corresponds to the values 
$(\kin)=(0, 0, 1, \dots, 1 )$;
we recall that in all these examples $\k_1=0$. 
The geometrical meaning of the
contractions $\Gc$, $m=1, \dots, N$, is to describe the behaviour of the space
${\cal S}
\equiv SO_{\kin} (N+1) / SO_{\kiin}(N)$ around a point, a line,
\dots, an hyperplane. In particular, within the inhomogeneous CK family
$\k_1 = 0$, only those contractions $\Gc, m=2, \dots, N$ may produce a 
different algebra. In other words, these inhomogeneous algebras can be
thought of as the result of a  `local' contraction (around
a point, $m=1$) 
which made to vanish the 
associated curvature, although they   
can still be contracted to describe the behaviour of the
space around a line, \dots, an hyperplane, and hence the remaining 
contractions $\Gc, m=2, \dots, N$ may be relevant. 
For instance, the non-relativistic limit, where the behaviour of spacetime 
geometry is approximated in the neighbourhood of a given (time-like) line
corresponds to the contraction where $\k_2 \to 0$. 

\medskip

A second order central element for the algebra
$\is_{\kiin} (N)$, coming after a specialisation to this case of a
suitable rescaling of the general CK Killing's form, reads 
\beq
{\cal{\extr C}}= \sum_{i=1}^{N-1}\k_{iN}\P_i^2  + \P_N^2\quad;
\label{ad}
\eeq
notice that this Casimir only involves generators from the abelian translation
subalgebra. 

\medskip

Summarising, we see that the graded contraction language allows us to describe
contractions simply by setting some parameters equal  to zero.
These contractions may still be described by the standard IW framework, 
although the graded contraction scheme is more economical and permits a 
unified discussion of the different contractions.

\subsect{Contractions and dimensional analysis}

The minimal possible approach to study
the dimensional structure in CK algebras is by enforcing the dimensional
homogeneity of the commutation relations in all algebras in the CK family. In
this approach, all generators as well as the structure constants are
dimensional, in such a way that these dimensions are the same in all CK
algebras. Consider the redefinition
$\J_{ab}=\eta_{ab}\J'_{ab}$ for any CK algebra. If we now want the second 
commutator in (\ref{aa}) to be preserved (we still choose the structure
constants equal to $1$ as dimensionless, for we are interested here in
algebras in the CK family, and not beyond), we need 
$\eta_{ab}\eta_{bc}=\eta_{ac}$, so that we see that $\eta_{ab}$ may be 
expressed as  $\eta_{ab}=\eta_{a+1}\eta_{a+2}\ldots\eta_{b}$. If we make now
this change in the first and third commutators, we get 
\beq
[\J'_{ab},\J'_{ac}]=\frac{\k_{ab}}{\eta_{ab}^2}\J'_{bc},\qquad
[\J'_{ac},\J'_{bc}]=\frac{\k_{bc}}{\eta_{bc}^2}\J'_{ab}. \label{cambio}
\eeq
In the special case when all $\k_i$ are different from
zero (the case of simple algebras), the choice $\eta_{ab}^2=|\k_{ab}|$
leads to the standard commutators for the $\J'_{ab}$ of the real form 
$so(p,q)$ of the specific algebra considered, with all non-zero structure 
constants equal to $\pm 1$, and the $\J'_{ab}$'s are dimensionless. 

In the general case (for generic CK algebras) by virtue of the above
redefinition, the generators in (\ref{aa}) have as dimensions
$[\J_{ab}]=[\k_{ab}]^{1/2}= \prod_{i=a+1}^{b} [\k_{i}]^{1/2}$. In this
approach, each constant $\k_a$ has dimensions, and if the dimension of the
generator $\J_{a-1 a}$ is written as $D_{a}^{-1}$, then it is clear that the
dimensions of each $\k_a$ are $[\k_a]=D_{a}^{-2}$, irrespective of $\k_a$ 
being zero or not. The dimension of each generator $\J_{ab}$
includes a factor $[\k_a]$ for each of the
$\k_a$-rectangles in (\ref{triangulo}) to which $\J_{ab}$ belongs.

\bigskip

Another possibility is
to allocate dimensions to generators and/or canonical parameters for each CK
algebra independently, in such a way as to make dimensionless all non-zero
structure constants in the algebra.
The idea of basing the dimensional analysis of a theory on the structure of its 
underlying Lie group/algebra has its 
roots in the well known examples of the 
Poincar\'e and Galilei groups, which are obtained by contracting with respect 
two dimensionful parameters, the de Sitter radius and the velocity of light $c$, 
and has been discussed in \cite{CSAnalDim}; see also \cite{AldAzb}. 

If for a simple Lie algebra in the CK
family (\ref{aa}) with non-zero
$\k_i$ constants we adopt this hypothesis, then as a consequence all the
generators of the algebra, as well as their associated canonical parameters 
are also without dimensions (as the $\J'_{ab}$'s in the first approach). 
If the same requirement is applied to a non-simple CK Lie algebra, then we get 
the result that some generators are also dimensionless, while others get a 
dimension.
For example, if this is done on a CK algebra with a {\em a single} $\k_{a}$ is
equal to zero, it is clear from the commutation relations that those
generators which acquire in this case a non-trivial dimension $D_a^{-1}$
($D_a$ is then the dimension of corresponding canonical parameter) are exactly
those inside the $\k_a$-rectangle corresponding to the
constant $\k_a$ which vanished in the triangular arrangement of generators. 

If there are two constants equal to zero, say
$\k_{a}=\k_{b}=0$, there will be two non-trivial dimensions, and so on.
Remark that now each $\k_a$ has a dimension which is still $D_a^{-1}$ when 
$\k_a=0$ but it is dimensionless when $\k_a\neq0$. In this alternative choice
the dimensions of the generators will still be given by  
$[\J_{ab}]=[\k_{ab}]^{1/2}= \prod_{i=a+1}^{b} [\k_{i}]^{1/2}$ but those
factors where $\k_i \neq 0$ are $[\k_{i}]=1$.  This situation is  exemplified
in the transition from the Poincar\'e to the  Galilei algebras.  These are
given by the values
$(\k_1,\k_2,\k_3,\k_4)=  (0,-1/c^2,1,1)$ and
$(0,0,1,1)$ respectively, and the arrangement of the Galilei and or
Poincar\'e generators written in the usual physical notation is given by: 

\noindent 
\beq
\begin{tabular}{cccc}
$H $ & $ P_1 $ & $  P_2  $ & $  P_3  $  \\ 
$  $ & $ K_1 $ & $  K_2  $ & $  K_3  $  \\ 
$  $ & $     $ & $ J_{12}$ & $ J_{13}$  \\ 
$  $ & $     $ & $       $ & $ J_{23}$  \\ 
\end{tabular}
\label{physNot}  
\quad.
\eeq
The assignment of dimensions made this way for Poincar\'e and  for  Galilei
algebras is
$$
\noindent
\begin{tabular}{cccc}
$ D_1^{-1}$ & $ D_1^{-1}$ & $ D_1^{-1}$ &  $ D_1^{-1}$ \\ 
$  $ & $ 1 $ & $ 1   $ & $  1   $  \\ 
$  $ & $     $ & $   1  $ & $    1  $  \\ 
$  $ & $     $ & $       $ & $    1  $  \\ 
\end{tabular}
\qquad\qquad 
\begin{tabular}{cccc}
$ D_1^{-1}$ & $ D_1^{-1}D_2^{-1} $ & $ D_1^{-1}D_2^{-1} $ & 
$ D_1^{-1}D_2^{-1} $ \\ 
$  $ & $ \phantom{D_1^{-1}}D_2^{-1}  $ & 
$ \phantom{D_1^{-1}}D_2^{-1}   $ & $ 
\phantom{D_1^{-1}}D_2^{-1}   $  \\ 
$  $ & $     $ & $   1  $ & $    1  $  \\ 
$  $ & $     $ & $       $ & $    1  $  \\ 
\end{tabular}
$$
which gives the single `length' dimension in
relativistic physics and the customary
$T, L$ dimensions of  non-relativistic physics ($D_1 \equiv T$ and $D_2 \equiv
LT^{-1}$, so 
$D_1 D_2 = L$). 

\bigskip

The relation between both perspectives to the dimensional analysis of
CK algebras is as follows: All constants $\k_a$ can be considered at the
beginning as dimensionful, and then all generators are also dimensionful.
However, when a given $\k_a \neq 0$, the dimension $D_a$ can be
removed by taking $\k_a$ as a pure number, which can be set equal to 
$\pm 1$; this  is tantamount to fixing the scale of the  generators or, in
other words, to measuring the associated group parameter in  terms of the
corresponding unit much in the same way as in a relativistic  theory we may
adopt units in which $c=1$ (\ie, $\k_2=-1$ above). In the former example,
setting $\k_3=1=\k_4$ may be understood as having hidden universal constants
in the theory (cf. \cite{LeviLeblond}).  However, once a dimensionful $\k_a$
has been set equal to zero (\ie, a  contraction has been made), the
generators in the corresponding box retain a  dimension $[\k_a]^{1/2}$ since
they cannot be re-scaled any longer. This is why some generators
in the former Galilei example retain the non-removable dimensions $D_1, D_2$,
while $D_2$ disappears in the Poincar\'e case while $D_3, D_4$ have already
disappeared in both cases\footnote{The above is not the only group-theoretical 
mechanism for the introduction of 
dimensions.
Where centrally extended groups are physically relevant, the dimensions of the 
two-cocycle realizing the extension play a role. 
For instance, in the 
$(1+1)$-dimensional extended Galilei group we find two 
parameters characterising the two-dimensional cohomology space, which 
correspond to the mass and a (constant) force.}.

\sect{Deformed N--dimensional affine CK algebras}

All the family of affine $N$--dimensional ($N\geq 2$) CK algebras
$\is_{\kiin} (N)$ can be endowed with a standard 
deformed Hopf algebra structure which has been called `quantum'
inhomogeneous CK structure and which has been given in
\cite{BHOSab,BHOScd}.  
In order to avoid repeating statements on
the index ranges, we will conform in Sections III and IV
to the following convention: the range of a latin index $i, j, k$ will
be $1, \dots, N-1$, and the index $N$ will be dealt with separately,
unless otherwise stated explicitly. Also, when two indices $i,j$
appear  in a generator, we will always assume that $i<j$.  

Let ${\cal A}$ be the algebra of the formal power series in the deformation 
parameter $\pardeform$ with coefficients in the enveloping algebra 
$\U(\is_{\kiin} (N))$ of the Lie algebra 
$\is_{\kiin} (N)$ of (\ref{ab}). 
Then the
coproduct, counit, antipode and deformed commutation relations of the
algebra 
$\U_\pardeform(\is_{\kiin}(N))$, which is a Hopf algebra, are given
by

\medskip

\noindent
1) {\it Coproduct}: \hfill
\beq\ba{l}
\co(\P_N)=1\otimes \P_N+ \P_N\otimes 1,
   \qquad \qquad 
\co(\P_i)=e^{-\frac \pardeform 2  \P_N}\otimes \P_i +
            \P_i\otimes e^{\frac \pardeform 2  \P_N},\\[0.3cm]
\co(\J_{ij})=1\otimes \J_{ij}+ \J_{ij}\otimes 1      \\[0.3cm]
\co(\J_{iN})=e^{-\frac \pardeform 2  \P_N}\otimes \J_{iN} +
    \J_{iN} \otimes e^{\frac \pardeform 2  \P_N} \!-
   \frac \pardeform 2  \sum_{s=1}^{i-1} \J_{si}e^{-\frac \pardeform 2  \P_N}
                                     \otimes \k_{iN} \P_s +
     \frac \pardeform 2 \sum_{s=1}^{i-1}\k_{iN} \P_s 
              \otimes e^{\frac \pardeform 2  \P_N} \J_{si}  \\[0.3cm]
\qquad\qquad\ 
  + \frac \pardeform 2 \sum_{s=i+1}^{N-1}\J_{is}e^{-\frac \pardeform 2  \P_N}
                          \otimes \k_{sN} \P_s
   -\frac \pardeform 2 \sum_{s=i+1}^{N-1}\k_{sN} \P_s \otimes 
                e^{\frac \pardeform 2  \P_N} \J_{is}\quad;
\label{ba}
\ea
\eeq

\noindent
2) {\it Counit}: 
\beq
\v(\P_i)=
\v(\P_N)=
\v(\J_{ij})=
\v(\J_{iN})=0 \quad;
\label{bb}
\eeq

\noindent
3) {\it Antipode}:
\beq\ba{l}
\gamma(\P_i) = -\P_i, \quad 
\gamma(\P_N) = -\P_N, \quad
\gamma(\J_{ij}) = -\J_{ij}, \quad 
\gamma(\J_{iN}) = -\J_{iN}-\k_{iN}(N-1)\frac {\pardeform}2\P_i \quad;
\label{bc}
\ea\eeq
(it may be written in a compact way as 
$\gamma(\X)=-e^{(N-1)\frac {\pardeform}2 \P_N} \X \, e^{-(N-1)\frac {\pardeform}2 \P_N}$),

\medskip

\noindent
4) {\it Deformed commutators}:
\beq\ba{l}
[\J_{iN},\P_{j}]=\delta_{ij}\frac 1\pardeform\sinh(\pardeform\P_N), \cr
[\J_{iN},\J_{jN}]=\k_{jN}\biggl\{\J_{ij}\cosh(\pardeform\P_N)
+\frac{\pardeform^2}4\biggl(\sum_{s=1}^{i-1}\k_{iN}\P_s\W_{sij}\cr
\qquad\qquad\qquad 
-\sum_{s=i+1}^{j-1}\k_{sN}\P_s\W_{isj}+\sum_{s=j+1}^{N-1}
\k_{sN}\P_s\W_{ijs}\biggr)    \biggr\},\quad i<j, 
\label{bd}
\ea\eeq
where
\beq
\W_{ijk}=\k_{ij}\P_i \J_{jk}-\P_j\J_{ik}+\P_k\J_{ij},\quad i<j<k,
\qquad i,j,k=1,\dots,N-1 \ \ \  \quad.
\label{lubanski}
\eeq
The remaining commutators are non-deformed and as given in (\ref{ab}). 
It may be checked that $(A,\ \co, \ \v,\  \gamma)$ satisfies the 
Hopf algebra axioms 
and hence eqs. (\ref{ba})-(\ref{lubanski}) may be taken as the definition of 
the deformation $\U_\pardeform(\is_{\kiin} (N))$ of 
$\U(\is_{\kiin} (N))$.
The parameter $\pardeform$ has a dimension inverse to that of $P_N$ 
so that the product $\pardeform \P_N$ is dimensionless, and may be
interpreted as the  parameter left after contracting 
the deformed Hopf algebra ${\cal U}_q({so}_{\k_1\ldots\k_N}(N+1))$  
by previously redefining $q$ in terms of $\pardeform$ and the contraction 
parameter.
However, the expression of the deformation 
${\cal U}_q({so}_{\k_1\ldots\k_N}(N+1))$  
in the `physical' basis is not known and this precludes us for the moment from 
deriving  (\ref{ba})-(\ref{lubanski}) by contracting its deformed simple 
parent algebra ${\cal U}_q({so}_{\k_1\ldots\k_N}(N+1))$.
Nevertheless, it may be seen that the deformed  Hopf algebra
$\U_\pardeform(\is_{\kiin} (N))$ is a  quantisation of the coboundary Lie
bialgebra 
$(\U(\is_{\kiin} (N)),r)$ generated by the (non-degenerate) classical 
$r$--matrix
\beq
r=\pardeform\sum_{s=1}^{N-1} \J_{sN}\wedge \P_s\quad.\label{bf}
\eeq
Due to the structure of $r$ and our convention about dimensions, it
turns out that $r$ is dimensionless, regardless of the values of the constants
$\k_i$ 
since the product of $J_{sN}$ and $P_s$ will always have the same dimensions 
of $P_N$.

We remark that the above deformation (\ref{ba})-(\ref{lubanski}) is not the 
only one possible for the $\is_{\kiin}(N)$ family.
However, it is distinguished by the fact that all its members present
deformed algebra and coalgebra sectors.

The quantum analogue of the second order Casimir (\ref{ad}) is expressed by
\beq
{\cal {\extr C}}_\pardeform=
\sum_{i=1}^{N-1} \k_{iN}\P_i^2 + 
\frac 4{\pardeform^2}\left[\sinh(\frac \pardeform 2  \P_N)\right]^2\quad.
\label{ch}
\eeq

As far as their action on the algebra generators is concerned, 
the quantum versions
$\Gq$ of the classical  IW contractions $\Gc$ are defined to coincide
with the classical one  (\ref{ae}). In particular, the generator
$\P_N$ is re-scaled by the corresponding
contraction parameter $\epsilon$ in any of the contractions in the family
$\Gq$, $\Gc(\P_N)\equiv\P_N'=\epsilon \P_N$ (for $\epsilon\to 0$).
This means that since one has to replace $\P_N$ in (\ref{ba}) by 
$\P'_N/\epsilon$, the exponents there will diverge. 
It is therefore natural to replace simultaneously $\pardeform$ by 
$\epsilon\pardeform'$, \ie\ to re-scale the deformation parameter 
by $\Gc(\pardeform)\equiv\pardeform'=\pardeform/\epsilon$ 
(all primes are removed after taking the contraction limit),
as the simplest possibility to preserve the coproduct
(\ref{ba}). 
Summing up, the quantum contraction $\Gamma^{(M)}_\pardeform$ is defined as
the result of taking the limit $\epsilon  \to 0$ in (\ref{ba})-(\ref{bd}) once
the transformations
\beq
\Gq(\epsilon , \X)=\Gc (\epsilon , \X)\quad,\quad
\X=(\P,\J)\quad,\quad
\Gq(\epsilon , \pardeform)= \pardeform/\epsilon 
\label{bh}
\eeq
are performed.
We conclude this section with three observations.  
First, we have made constant reference to the IW procedure only because up to 
now the graded contraction theory had not been extended to deformed algebras. 
Second, as far as the generators $\X$ are concerned, 
$\Gq=\Gc$, so that only the action of 
$\Gq$ on $\pardeform$ makes  $\Gc$ and $\Gq$ different.
The third comment is that the re-scaling of $\pardeform$ and $\P_N$
implied by  $\Gq$ may change their dimensions (see Sec. 2) while consistently
keeping a dimensionless $\pardeform\P_N$ exponent.


\sect{Bicrossproduct structure of $\U_\pardeform(\is_{\kiin} (N))$}

It is not obvious to see 
whether the Hopf algebra $\U_\pardeform(\is_{\kiin}(N))$ 
has a bicrossproduct structure by a simple inspection of
(\ref{ba})--(\ref{bd}). The clue in this direction is provided by the
bicrossproduct structure \cite{KPoinBicros} of the
$\kappa$-Poincar\'e algebra \cite{LNRT} 
(appearing in our scheme when $(\k_1,\k_2, \k_3, \k_4) = (0, 1, 1, -1)$), 
which is clearly displayed in terms of a new set of
generators.

The aim of this section is to show that {\it all} the deformed
Hopf algebras  
$\U_\pardeform(\is_{\kiin} (N))$ in the CK family have indeed a 
bicrossproduct structure. The basic bicrossproduct 
formulae used are recalled in the Appendix; for a
detailed exposition, see \cite{MajMajid}.
Let ${P}_i$, ${J}_{ij}$, $P_N$, $J_{iN}$
be the new set of generators defined in terms of the old ones ${\P}_i$ and 
${\J}_{ij}$ by 
\beq
\ba{l}
{P}_{i}= e^{-(\pardeform/2)\P_N}\P_{i}, 
\qquad {P}_{N}=\P_{N}, \qquad  
{J}_{ij}=\J_{ij},\\[0.3cm]
{J}_{iN}=\frac{1}{2}\{\J_{iN},e^{-(\pardeform/2)\P_N}\} + 
  \frac{\pardeform}{4}\sum_{s=1}^{i-1} \k_{iN}\{\J_{si},\P_s\}e^{-(\pardeform/2)\P_N} - 
  \frac{\pardeform}{4}\sum_{s=i+1}^{N-1}\k_{sN}\{\J_{is},\P_s\}e^{-(\pardeform/2)\P_N}\quad, 
\label{ca}
\ea\eeq
A straightforward but tedious computation leads to the following 
new expressions 
(where $i,j,k=1,\ldots,N-1$) for the coproduct (\ref{ba}), counit (\ref{bb}), 
antipode (\ref{bc}) and algebra commutators ((\ref{bd}) and/or (\ref{ab})): 

\noindent
1) {\it Coproduct}: \hfill
\beq\ba{l}
\co(P_i)=e^{-\pardeform P_N}\otimes P_i + P_i\otimes  1 \quad, \quad  
\co({P_N})=1\otimes{P_N}+{P_N}\otimes 1 \quad,
\\[0.3cm]    
\co({ J_{ij}})=1\otimes{ J_{ij}}+{J_{ij}}\otimes 1\quad,\\[0.3cm]
\co(J_{iN})=e^{-\pardeform P_N}\otimes J_{iN} +J_{iN}\otimes 1 
  +\pardeform\sum_{s=1}^{i-1}\k_{iN} P_s \otimes J_{si}
  -\pardeform\sum_{s=i+1}^{N-1}\k_{sN} P_s \otimes J_{is},
\label{cb}
\ea\eeq

\noindent
2) {\it Counit}: 
\beq
\v(P_i)=\v(P_N)=\v(J_{ij})=\v(J_{iN})=0,
\label{cc}
\eeq

\noindent
3) {\it Antipode}:
\beq\ba{l}
\gamma(P_i)=-\exp(\pardeform P_N) P_i\quad,\quad 
\gamma(P_N)=-P_N \quad,\quad 
\gamma(J_{ij})=-J_{ij}\quad,
  \\[0.3cm]
\gamma(J_{iN}) = - e^{\pardeform P_N}J_{iN} + 
 \pardeform e^{\pardeform P_N}\sum_{s=1}^{i-1}\k_{iN} P_s J_{si} -
 \pardeform e^{\pardeform P_N}\sum_{s=i+1}^{N-1}\k_{sN} P_s J_{is}, 
\label{cd}
\ea\eeq

\noindent
4) {\it Commutators}: 
\begin{eqnarray}
&& [P_i,P_j]=0 \quad, \qquad [P_i,P_N]=0\quad, \nonumber \\{} 
&& [J_{ij},J_{ik}] =  \k_{ij} J_{jk} \qquad,
   [J_{ij},J_{jk}] = -J_{ik}   \quad,
   [J_{ik},J_{jk}] = \k_{jk}J_{ij}, \nonumber \\{} 
&& [J_{ij},J_{iN}] =  \k_{ij} J_{jN}, \qquad
   [J_{ij},J_{jN}] = -J_{iN} ,  \qquad
   [J_{ik},J_{jN}] = \k_{jN}J_{ij}, \label{ce} \\{} 
&& [J_{ij},P_k]=\delta_{ik}P_k - \delta_{jk}\k_{ij}P_i, \qquad
   [J_{ij},P_N]=0, \nonumber \\{} 
&& [J_{iN},P_{j}]=\delta_{ij}\left(\frac{1-e^{-2\pardeform P_N}}{2\pardeform}- 
     \frac \pardeform 2  \sum_{s=1}^{N-1}\k_{sN}P_s^2 \right) 
     + \pardeform \k_{iN}P_iP_j\quad, \qquad  
   [J_{iN},P_N]=-\k_{iN}P_i \quad . \nonumber 
\end{eqnarray}
Thus, all brackets for the new generators $ P_i, P_N, J_{ij}, J_{iN} $ 
coincide with the non-deformed ones given in (\ref{ab}) 
(substituting  everywhere the new $X$'s for their counterparts $\extr X$) 
except for $[J_{iN},P_{j}]$,  
which is now the only deformed commutation relation.
The effect of (\ref{ca}) is to modify the second commutator in 
(\ref{bd}), so that one recovers the undeformed $\s_{\kin}(N+1)$
algebra commutators, and to replace the commutators in the first line in
(\ref{bd}) by those in the last line in (\ref{ce}).  As a result, terms with
the
$\W$ symbols are no longer present in the deformed commutators.

It may be checked that for  $\pardeform={1\over\kappa}$ and $N=4$ with 
$(\k_1,\k_2,\k_3,\k_4)=(0,1,1,-1)$ eqs. (\ref{cb})-(\ref{ce}) reproduce the 
$\kappa$-Poincar\'e algebra in the basis of \cite{KPoinBicros} 
for which $[\kappa]=L^{-1}$, $[P_4]=L^{-1}$.
If we want $P_4$ to have dimensions of inverse of time we may take 
$(\k_1,\k_2,\k_3,\k_4)=(0,1,1,-c^2)$ instead
since $\k_1$, before being set equal to zero, was 
$\k_1={1\over R^2}$; in this case $[\kappa]=T^{-1}$. 
We check that the metric after (\ref{ac}) will diverge in a non-relativistic 
limit with $\k_4=-c^2$, which explains why a non-relativistic limit of the 
$\kappa$-Poincar\'e algebra \cite{LNRT} requires a further redefinition of the 
deformation parameter $\kappa$ (see the end of Sec. 5).

The new expressions for the coproduct, counit, antipode and 
commutation relations of $\U_\pardeform(\is_{\kiin} (N))$ now allow us to 
uncover its bicrossproduct structure.
To this aim, consider
the translation sector, generated by $\{ P_1, \dots, P_N \}$.
According to the expressions (\ref{cb})--(\ref{ce}),  it defines a commutative
but non co-commutative Hopf subalgebra of 
$\U_\pardeform(\is_{\kiin} (N))$ which will be denoted as
$\U_\pardeform(T_N)$.  Let now $\U( {so}_{\kiin} (N))$ be the 
non-commutative and co-commutative non-deformed  
CK Hopf algebra spanned by the remaining
generators 
$\{J_{ij};\ i<j,\  i,j=1,\dots, N\}$, hence with commutation relations given by 
(\ref{ab}) and primitive coproduct (when all
$\k$'s are non-zero, this is a pseudo-orthogonal algebra). Let us define a
right action 
$\a:\U_\pardeform(T_N)\otimes \U({so}_{\kiin} (N)) \to \U_\pardeform(T_N)$ by 
\beq
\a(P_i,J_{jk})\equiv P_i \lhd J_{jk} :=[P_i,J_{jk}] ; \quad
j<k,\ \ i,j,k= 1,2,\dots, N,\label{cf}
\eeq       
where the commutators is given in (\ref{ce}), 
and a left coaction $\b:\U({so}_{\kiin} (N)) \to \U_\pardeform(T_N)\otimes 
\U({so}_{\kiin} (N))$ by
\beq\ba{l}
\b(J_{ij}) :=1\otimes J_{ij}\quad,\quad\\[0.3cm]
\b(J_{iN}) :=e^{-\pardeform P_N}\otimes J_{iN} +
\pardeform\sum_{s=1}^{i-1}\k_{iN} P_s \otimes J_{si}
-\pardeform\sum_{s=i+1}^{N-1}\k_{sN} P_s \otimes J_{is}\quad.
\ea \label{cg} \eeq
It may be checked that $\U_\pardeform(T_N)$ is a right 
$\U({so}_{\kiin} (N))$--module 
algebra ($\U({so}_{\kiin} (N))\rimo\U_\pardeform(T_N)$) and 
that $\U({so}_{\kiin}(N))$ is 
a left $\U_\pardeform(T_N)$-comodule coalgebra  ($\U({so}_{\kiin} 
(N))\leco\U_\pardeform(T_N)$) under the action (\ref{cf}) and coaction (\ref{cg}),
respectively,
and that the compatibility conditions \cite{MajMajid} 
(\ref{apix})--(\ref{apxiii}) between $\alpha$ y 
$\beta$ needed for
$\U({so}_{\kiin} (N))\otimes \U_\pardeform(T_N)$ to have a
bicrossproduct  structure are fulfilled.
For instance, (\ref{apxiii}) is automatically satisfied, since
$\U(\s_{\kiin}(N))$ is undeformed and hence co-commutative and $\U_\pardeform(T_N)$ 
is abelian.
This case, specially relevant here, was discussed in \cite{SingMoln}.
Then, the bicrossproduct structure of $\U_\pardeform(\is_{\kiin} (N))$ may
be stated  in the form of the following 

\bigskip
\noindent
{\bf Theorem}. 
{\it
The deformed Hopf CK family of algebras $\U_\pardeform(\is_{\kiin} (N))$ has a
bicrossproduct structure 
$$
\U_\pardeform(\is_{\kiin} (N)) =  \U({so}_{\kiin} (N))
^\beta\bicross_\alpha  
\U_\pardeform(T_N)
$$
relative to the right action $\alpha$ and left coaction $\beta$ given by 
(\ref{cf}) and (\ref{cg}) respectively.}

\medskip
\noindent
{\sl Proof}:
As mentioned, the mappings $\a$ and $\b$ satisfy the bicrossproduct 
conditions as may be checked by direct computation. 
Then the expressions
(\ref{apxiv})-(\ref{apxvii})  give the associated coproduct, counit and 
antipode. 
It is then verified that the resulting expressions 
are in agreement with (\ref{cb})--(\ref{ce}) {\it q.e.d.}

The interesting consequence of the above discussion is that, as the direct 
inspection of expressions (\ref{cf}) and
(\ref{cg}) shows, the action and the coaction mappings depend on
the parameters $\k_i$ in such a way that the bicrossproduct structure
is formally invariant under any contraction $\k_i=0$. 
In other words, for $\pardeform$-deformations in the affine CK family,
the bicrossproduct structure is preserved by {\it all} 
the successive contractions: 
contracting and taking bicrossproduct of the appropriate Hopf algebras with the 
resulting actions and coactions are commuting processes. 
This is well within the
spirit of the CK scheme, the aim of which is to state properties which
hold simultaneously for a large number of algebras.

The expression of the deformed Casimir (\ref{ch}) in the new basis 
is
\beq
{\cal C}_\pardeform=
\sum_{i=1}^{N-1} \k_{iN}e^{-\pardeform P_N}{P}_i^2 + 
\frac 4{\pardeform^2}\left[\sinh(\frac \pardeform 2  {P}_N)\right]^2
\quad;
\eeq
it only depends on the generators of the 
deformed Hopf subalgebra $\U_\pardeform(\T_N)$.

On the other hand,  the expression for the $r$--matrix is similar to
the former (\ref{bf}) but in terms of the new generators, 
\beq
r=\pardeform\sum_{s=1}^{N-1} J_{sN}\wedge P_s\quad. \label{bfnew}
\eeq

\sect{Applications}
The quantum algebras we are dealing with range from deformations of the
inhomogeneous algebras 
${iso}(p,q), \ p+q=N$ (when all constants $\k_2, \k_3, \dots \k_N$
are different from zero) to the extreme case of a Hopf deformation of 
the algebra, where all constants are equal to zero, which
can be called flag space algebra 
(in this case the group action preserves a complete flag). 

Classically, all these algebras are semidirect products, and indeed there
is a semidirect structure in the CK algebras associated to the
vanishing of each constant $\k_i$. We have restricted here to the
algebras with $\k_1=0$, all of which have the semidirect structure
displayed in (\ref{ab}). When {\em all} remaining constants
$\k_i$ are different from zero, say $\k_i=\pm 1$, 
the algebra $\is_{\kiin}(N)$ is
isomorphic to an inhomogeneous pseudo-orthogonal algebra
${iso}(p,q), \  p+q=N$, 
with the semidirect structure given by  the natural
action of
${so}(p,q)$ on $\R^N$. These algebras are physically very relevant
and some of their deformations have been thoroughly studied. 
In particular, the $\pardeform$-deformed structures given in
(\ref{ba})--(\ref{bd}) include a deformed 
$N$-dimensional Euclidean algebra, a deformed $(N-1,1)$ Galilei algebra and
several
deformed $(N-1,1)$ Poincar\'e algebras, as well as their analogues for any
signature. 

The action and coaction mappings associated with the bicrossproduct 
are given by (\ref{cf}) and (\ref{cg}).
Explicitly, eq. (\ref{ce}) gives
\begin{eqnarray}
&& \a(P_N,J_{ij})\equiv P_N \lhd J_{ij} := 0 \quad,\quad
   \a(P_N,J_{iN})\equiv P_N \lhd J_{iN} := \k_{iN} P_i \quad,\nonumber \\
&& \a(P_k,J_{ij})\equiv P_k \lhd J_{ij} := 
            -\delta_{ki} P_j + \delta_{kj} \k_{ij} P_i\quad,  \nonumber \\
&& \a(P_k,J_{iN})\equiv P_k \lhd J_{iN} :=  
      -\delta_{ki}\left(\frac{1-e^{-2\pardeform P_N}}{2\pardeform}- \frac \pardeform 2  
       \sum_{s=1}^{N-1}\k_{sN}P_s^2 \right) - \pardeform \k_{iN}P_iP_k\quad.
        \label{includ}
\end{eqnarray}

If we consider the special case $N=4$, this set of algebras
includes {\em four} deformed Poincar\'e algebras 
$\U_\pardeform({\fra p}^{(s)}(3,1)), \  s=1,2,3,4$. These are
deformations of the four undeformed CK algebras, denoted ${\fra p}^{(s)}(3,1),\ 
s=1,2,3,4$, which are isomorphic to the 
(3,1) Poincar\'e algebra, and correspond to
identifying one of the generators $P_i$ to the time translation generator, the other
three being space translations.  
If the time generator is taken successively to be our  
$P_1$, $P_2$, $P_3$, $P_4$, these four algebras correspond to the four sets
of values of  $(\k_1,\k_2,\k_3,\k_4)=$ 
$(0,-1/c^2,1,1)$,
$(0,-c^2,-1/c^2,1)$, 
$(0,1,-c^2,-1/c^2)$, and
$(0,1,1,-c^2)$. The set of four deformed Poincar\'e algebras 
$\U_\pardeform({\fra p}^{(s)}(3,1)), \  s=1,2,3,4$ 
\cite{BHOScd} includes three `space-like' Poincar\'e deformed algebras, 
the last one being the $\kappa$-Poincar\'e algebra once  
$\pardeform={1/\kappa}$ with
$[P_N]=T^{-1}$). In each case, the rotation generators comprise the
boost and space rotation generators and the identification is made
according to the choice of the time generator (\eg, in the 
$\kappa$-Poincar\'e the boosts are the $N_i=J_{i4}$).  
The $N$-dimensional $\kappa$--Poincar\'e  \cite{MasLuk} 
is associated to the $\k_i$ values  $(0,1,\dots, 1,-c^2)$.

The Euclidean algebra $\fra e(4)$ appears only once (up to rescalings) for 
$(\k_1,\k_2,\k_3,\k_4)=(0,1,1,1)$
and the bicrossproduct structure of their Hopf CK quantum deformation  
$\U_\pardeform(\fra e(4))$, is  
$$
\U_\pardeform(\fra e(4))=  \U({so}(4)) \bicross \U_\pardeform(\R_4)
$$
(see \cite{APb} in the lower dimensional case).

The remaining quantum Hopf algebras in the CK family are quantum
deformations of less known undeformed Lie algebras, yet their
bicrossproduct structure is described in parallel to the former cases.
Specially relevant from the physical point of view 
is the Galilei algebra 
$\fra g (3,1)$, appearing within the affine CK family for the $\k$ values
$(0,0,1,1)$ (and only for these). It is worth remarking that this
Galilei algebra is obtained from the 
$\U_\pardeform({\fra p}^{(1)} (3,1))$ 
associated to the $\k$ values $(0,-1/c^2,1,1)$ by means of the
contraction $\k_2
\to 0$ (\ie, $c \to \infty$), which gives a deformation different from that in 
\cite{POLACOS,APa}. 
The bicrossproduct structure of the resulting $\U_\pardeform(\fra g (3,1))$
Hopf algebra (obtained for $(0,0,1,1)$ is  
$$
\U_\pardeform(\fra g (3,1)) = 
\U({so}_{(0,1,1)} (4)) \bicross \U_\pardeform(\R_4)
\equiv 
\U({iso}_{(1,1)}(3))\bicross  \U_\pardeform(\R_4)\quad.
$$  

It is interesting to check how the action and coaction mappings for the
Poincar\'e Hopf algebra $\U_\pardeform (\fra p^{(1)} (3,1))$ 
reduce to the corresponding Galilean ones
under the contraction $\k_2 \to 0$: $\U_\pardeform (\fra p^{(1)} (3,1))$ to
$\U_\pardeform (\fra g (3,1))$. 
We give explicitly the complete Hopf structure of both 
$\U_\pardeform({\fra p}^{(1)}(3,1))$ and $\U_\pardeform({\fra g}(3,1))$, which
correspond to the choices $(0,-1/c^2,1,1)$ and $(0,0,1,1)$ for the $\k_i$'s. We
will present the results in the usual physical basis, constituted by the 
generators of time translation,
$H$, space translations $P_1, P_2, P_3$, boosts $K_1, K_2, K_3$, and
space rotations $J_1, J_2, J_3$, 
which are related to the CK original generators $\{ P_i,\ J_{ij};\
i,j=1,2,3,4\}$ as follows: the three {\em space} translations, now denoted
as $P_1, P_2, P_3$ in order to conform with the standard physical notation, 
correspond to those formerly denoted as $P_2, P_3, P_4$ in (\ref{cb})-(\ref{ce}),
while the time translation generator
$H$  now corresponds to the former $P_1$
and the rest in (\ref{cb})-(\ref{ce}) correspond to 
$J_{12}=K_1\,,\,J_{13}=K_2\,,\,J_{14}=K_3\,,\,
J_{34}=J_1\,,\,J_{24}=-J_2\,,\,J_{23}=J_3$, as given in the diagram
(\ref{physNot}).  
 
The Hopf structure of the Poincar\'e algebra 
${\cal U}_\pardeform(\fra p^{(1)} (3,1))$ now follows 
from expressions (\ref{cb})-(\ref{ce}):

\medskip
\noindent
1) {\it Coproduct}: \hfill
\beq\ba{l}
\co(X)=e^{-\pardeform P_3}\otimes X + X\otimes  1, \qquad X\in \{H,P_1,P_2\},
\\[0.3cm]    
\co({P_3})=1\otimes{P_3}+{P_3}\otimes 1, \\[0.3cm]
\co({X})=1\otimes{X}+{X}\otimes 1, \qquad X\in \{K_1,K_2,J_3\},
\\[0.3cm] \co(K_{3})=e^{-\pardeform P_3}\otimes K_{3} +K_{3}\otimes 1 
  -\pardeform  P_1 \otimes K_{1}
  -\pardeform P_2 \otimes K_{2}, \\[0.3cm]
\co(J_{1})=e^{-\pardeform P_3}\otimes J_{1} +J_{1}\otimes 1 
 +\pardeform  H \otimes K_{2}
  +\pardeform P_1 \otimes J_{3}, \\[0.3cm]
\co(J_{2})=e^{-\pardeform P_3}\otimes J_{2} +J_{2}\otimes 1 
 -\pardeform  H \otimes K_{1}
  +\pardeform P_2 \otimes J_{3}. 
\label{fb}
\ea\eeq

\noindent
2) {\it Counit}: 
\beq
\v(H)=\v(P_i)=\v(K_{i})=\v(J_{i})=0,\qquad i=1,2,3.
\label{fc}
\eeq

\noindent
3) {\it Antipode}:
\beq\ba{l}
\gamma(X)=-e^{\pardeform P_3} X,\qquad X\in \{H,P_1,P_2\}; \qquad 
\gamma(P_3)=-P_3, \\ [0.3cm]
\gamma(X)=-X,\quad X\in \{K_1,K_2,J_3\}, \\[0.3cm]
\gamma(K_{3}) = - e^{\pardeform P_3}K_{3}
- \pardeform e^{\pardeform P_3} P_1 K_{1} 
- \pardeform e^{\pardeform P_3}P_2 K_{2}, \\[0.3cm] 
\gamma(J_{1}) = - e^{\pardeform P_3}J_{1} 
+ \pardeform e^{\pardeform P_3} H K_{2} 
+ \pardeform e^{\pardeform P_3} P_1 J_{3}, \\[0.3cm] 
\gamma(J_{2}) = - e^{\pardeform P_3}J_{2} 
- \pardeform e^{\pardeform P_3} H K_{1} 
+ \pardeform e^{\pardeform P_3}P_2 J_{3}.\label{fd}
\ea\eeq

\noindent
4) {\it Commutators}: 
\beq
\begin{array}{l}
[H,P_i]=0,\quad [P_i,P_j]=0\quad i,j=0,1,2,3 ,\\[0.3cm]
[K_{i},H]=P_i,\ \ \ i=1,2; \qquad 
   [K_{3},H]=\frac{1-e^{-2\pardeform P_3}}{2\pardeform}  
   -\frac{\pardeform}{2c^2}H^2- \frac{\pardeform}{2}P_1^2 
   -\frac{\pardeform}{2}P_2^2, \\[0.3cm]
[K_{i},P_j]=\frac{1}{c^2} \delta_{ij}H, \quad 
i,j=1,2, \qquad [K_{i},P_3]=0,\\[0.3cm] 
[K_{3},P_1]=-\frac{\pardeform}{c^2}HP_1,\qquad  
  [K_{3},P_2]=-\frac{\pardeform}{c^2}HP_2,\qquad
  [K_{3},P_3]=\frac{1}{c^2}H, \\[0.3cm]
[J_{1},H]={\pardeform}P_2H,\qquad
  [J_{2},H]=-{\pardeform}P_1H,\qquad [J_{3},H]=0,\\[0.3cm]
 [J_{1},P_1]={\pardeform}P_2P_1, \quad 
  [J_{1},P_2]=\frac{1-e^{-2\pardeform P_3}}{2\pardeform}  
  +\frac{\pardeform}{2c^2}H^2- \frac{\pardeform}{2}P_1^2 
  +\frac{\pardeform}{2}P_2^2,\quad 
  [J_{1},P_3]=-P_2,  \\[0.3cm]
 [J_{2},P_1]=-\frac{1-e^{-2\pardeform P_3}}{2\pardeform}  
  -\frac{\pardeform}{2c^2}H^2 - \frac{\pardeform}{2}P_1^2 
  +\frac{\pardeform}{2}P_2^2, \quad  [J_{2},P_2]=-{\pardeform}P_1P_2,
  \quad [J_{2},P_3]=P_1, \\[0.3cm]
[J_{3},P_1]=P_2, \quad [J_{3},P_2]=-P_1, \quad[J_{3},P_3]=0, 
  \\[0.3cm]  
[K_{i},K_j]=-\frac{1}{c^2}\varepsilon_{ijk}K_k, \quad 
  [J_{i},K_j]=\varepsilon_{ijk}K_k,\quad
  [J_{i},J_j]=\varepsilon_{ijk}J_k,\ \ \ i,j=1,2,3. 
\end{array}
\label{fe}
\eeq

The explicit form of the action and the coaction for $\U_\pardeform({\fra
p}^{(1)}(3,1))$, given by (\ref{cf}) and (\ref{cg}) is obtained from  the
commutators (\ref{fe}) and the coproduct (\ref{fb}). The bicrossproduct structure
of   $\U_\pardeform({\fra p}^{(1)}(3,1))$ is
$$
\U_\pardeform({\fra p}^{(1)}(3,1)) = 
\U({so}_{(-1/c^2,1,1)} (4)) \bicross \U_\pardeform(\R_4).
$$  

In the non-relativistic limit, $\k_2=-1/c^2 \to 0$, this deformed algebra goes 
to a 
new deformed Galilei algebra, $\U_\pardeform({\fra g}(3,1))$, whose coproduct,
counit and antipode are the same as 
in  (\ref{fb})--(\ref{fd})). With respect to the Lie commutators we write
only those that are different from the Poincar\'e 
$\U_\pardeform({\fra p}^{(1)}(3,1))$ ones in (\ref{fe})

\noindent
\begin{eqnarray}
&&  [K_{3},H]=\frac{1-e^{-2\pardeform P_3}}{2\pardeform}  
- \frac{\pardeform}{2}P_1^2 
-\frac{\pardeform}{2}P_2^2,\nonumber\\[0.1cm]
&&[K_{i},P_j]=0, \quad 
i, j=1,2,3,\nonumber\\[0.1cm] 
&& [J_{1},P_2]=\frac{1-e^{-2\pardeform P_3}}{2\pardeform}  
- \frac{\pardeform}{2}P_1^2 
+\frac{\pardeform}{2}P_2^2, 
\label{ff}
\\[0.1cm]
&&[J_{2},P_1]=-\frac{1-e^{-2\pardeform P_3}}{2\pardeform}  
 - \frac{\pardeform}{2}P_1^2 
+\frac{\pardeform}{2}P_2^2,\nonumber\\[0.1cm]
&&[K_{i},K_j]=0,\quad  i,j=1,2,3. \nonumber
\end{eqnarray}
{}From simple inspection we see that the 
Galilean action (as the commutators) has
changed from 
that in $\U_\pardeform({\fra p}^{(1)}(3,1))$,
yet the coaction is the same than in
Poincar\'e since the coproduct has remained the same.

All other Poincar\'e quantum algebras (including the $\kappa$-Poincar\'e) do 
not allow a direct `non-relativistic' contraction; this is clear when the 
constant $c$ is explicitly written as above (see \cite{APa} for a 
discussion).

\sect{The group deformation aspect: the case of \hfil\break
      $\Fun_\pardeform (ISO_{\k_2}(2))$}

The bicrossproduct structure of 
$\U_\pardeform (\is_{\kiin}(N))$ opens the way 
to the possibility of recovering more easily the dual group 
$\Fun_\pardeform(ISO_{\kin} (N))$ from the dual bicrossproduct `group-like'
expressions. We show here the explicit calculation in the lowest
dimension $N=2$ to exhibit the procedure. 
As shown in Sec. 4, the deformed Hopf algebra has the bicrossproduct
structure 
\beq 
\U_\pardeform (\is_{\k_2}(2)) =  \U({so}_{\k_2} (2)) \bicross 
 \U_\pardeform (\T_2)\quad.
\label{da}
\eeq
Let us now recover from it the quantum dual group, 
$\Fun_\pardeform(ISO_{\k_2}(2))$. 

The dual algebra $\Fun_\pardeform(\T_2)$ of $\U_\pardeform(\T_2)$ is easily found from 
$\U_\pardeform(\T_2)$ to be
\beq
\co(a_1)=1\otimes  a_1 +  a_1\otimes 1\quad,\quad
\co(a_2)=1\otimes  a_2 +  a_2\otimes 1\quad; \label{db}
\eeq
\beq
\v(a_i)=0\quad;\quad \gamma(a_i)=-a_i\quad,\quad i=1,2\quad;\label{dc}
\eeq
and 
\beq
[a_1,a_2]=\pardeform a_1\quad,
\eeq
where $a_1,a_2$ are a system of coordinates of $\Fun_\pardeform (\T_2)$. 

On the other hand, let 
$\Fun(SO_{\k_2}(2))$ be the dual of
$\U(SO_{\k_2}(2))$ generated by $\varphi$,  
with a non-deformed Hopf structure defined by
\beq 
\co(\varphi)= 1\otimes \varphi +  \varphi\otimes 1\quad, \quad 
\v(\varphi)=0\quad,\quad \gamma(\varphi)=-\varphi\quad. \label{dd} 
\eeq   
Now, the problem is finding a pair of mappings $\bar{\b}$ and
$\bar{\a}$ (see Appendix),  
\beq\ba{l}
\bar{\b}\ :\ \Fun_\pardeform (\T_2)\to \Fun_\pardeform (\T_2)\otimes
\Fun(SO_{\k_2}  (2)),\\[0,3cm]
\bar{\a} \ :\ \Fun_\pardeform (\T_2)\otimes \Fun(SO_{\k_2}  (2))\to 
\Fun(SO_{\k_2} (2)).\label{de}
\ea\eeq 
duals, respectively, to $\a$ and $\b$ as given in (\ref{includ}), (\ref{cg}) 
for $N=2$. A calculation shows that they have the form 
\beq\ba{l}
\bar{\b}(a_1)=a_1\otimes \C_{\k_2}(\varphi) +
a_2\otimes {\k_2} \S_{\k_2}(\varphi)\quad,\\[0.3cm]
\bar{\b}( a_2)=-a_1\otimes \S_{\k_2}(\varphi) +
a_2\otimes  \C_{\k_2}(\varphi)\quad,\label{df}
\ea\eeq
and 
\beq
\bar{\a}( a_1\otimes \varphi)=\pardeform(1- \C_{\k_2}(\varphi))
\quad,\quad
\bar{\a}(a_2\otimes \varphi)=\pardeform\S_{\k_2}(\varphi)
\quad,
\label{dh}
\eeq
where the functions $\C_{\k}(\varphi)$ and $\S_{\k}(\varphi)$ reduce to
the trigonometric cosine and sine functions for $\k=1$ and to the
hyperbolic ones for
$\k=-1$ (see \cite{BHOSe} for more details). Note that 
$$
(\bar{\b}( a_1),\bar{\b}( a_2))=(a_1,a_2)\dot\otimes
\left( \begin{array}{rr} 
\C_{\k_2}(\varphi) & -\S_{\k_2}(\varphi) \\\k_2 \S_{\k_2}(\varphi) & 
\C_{\k_2}(\varphi) \end{array}  \right)\quad,\label{dj}
$$
where the $2\times 2$ matrix is the transpose of the matrix
representing the generic element of the group $SO_{\k_2} (2)$.   

Since $\bar\beta$ modifies the originally co-commutative coproduct in 
$\Fun_\pardeform(\T_2)$ and $\bar\alpha$ the commutation relations between the 
generators of the two algebras $\Fun_\pardeform(\T_2)$ and that generated by 
$\varphi$ (see eqs. (\ref{abpxv}) and (\ref{abpxiv})),
we may now determine through its bicrossproduct structure the deformed Hopf 
algebra $\Fun_\pardeform (SO_{\k_2} (2+1))$ which has the form 
\beq\ba{l}
\co( a_1)=1\otimes a_1 +
a_1\otimes  \C_{\k_2}(\varphi)+  a_2\otimes
{\k_2}\S_{\k_2}(\varphi),\\[0.3cm]
\co( a_2)=1\otimes a_2 -
a_1\otimes  \S_{\k_2}(\varphi) + a_2\otimes
\C_{\k_2}(\varphi),\\[0.3cm]
\co(\varphi)=1\otimes \varphi + \varphi \otimes  1;\label{dk}
\ea\eeq
\beq
\v(a_i)=0\quad,\ \ \  i=1,2\quad;\quad \v(\varphi)=0\quad;\label{dl}
\eeq
\beq\ba{l}
\gamma(a_1) = 
    - \C_{\k_2}(\varphi) a_1 -{\k_2} \S_{\k_2}(\varphi) a_2\quad,  \\[0.3cm]
\gamma(a_2) = 
     \S_{\k_2}(\varphi) a_1 -\C_{\k_2}(\varphi) a_2\quad;
\label{dm}
\ea\eeq
\beq
[a_1,\varphi]=  \,\pardeform \,(1-\C_{\k_2}(\varphi))\quad,\quad
[ a_2,\varphi]= \,  \pardeform \S_{\k_2}(\varphi)\quad,\quad
[{a_1},{a_2}]=  \,\pardeform  \,  a_1\quad.\label{dn}
\eeq
In this way the results obtained in \cite{BHOSe} (and \cite{APb} for $\k_2=1,\ 
a_1\to -a_2,\ a_2\to a_1$) are recovered.


\sect{Conclusions}

The bicrossproduct structure holds
true for all deformed Hopf algebras in the inhomogeneous CK family 
$\U_\pardeform (\is_{\kiin}(N))$. 
The theorem in Sec. 4 follows from the fact that 
the action and the coaction mappings that characterise the
bicrossproduct structure of $\U_\pardeform(iso_{\k_2\dots\k_N}(N))$ depend on 
$\k_2,\k_3,\ldots \k_N$ in such a way that the contractions are simply
described by setting (some of) them equal to zero. 
In this sense we may say that
the bicrossproduct structure  (\ref{cf})--(\ref{cg}) is compatible with
the contractions  (\ref{bh}).  As a result, any contracted
$\pardeform$-algebra is also  a bicrossproduct for  the `contracted'
action and coaction.  We can also state this result by saying that the
bicrossproduct structure does not diverge under any of the 
contractions within the CK affine family. 
These results extend obviously to the dual 
$\Fun_\pardeform(ISO_{\k_2,\ldots,\k_N}(N))$ family of deformed Hopf algebras.

The deformation parameter $\pardeform$ behaves in the same way under the 
complete family of CK contractions of the inhomogeneous algebras
$\U_\pardeform(\is_{\kiin} (N))$; if the (first, for instance) 
contraction has been realized as a limit,
the deformation parameter $z$ in $q=e^z$ is redefined by 
$\Gq(z)=z'={\lambda\over R}$ 
(cf. eq. (\ref{bh})) and the limit corresponds to $[\k_1]^{1/2}=1/ R\to 0$.
As discussed, this can be understood as a mechanism assigning 
dimensions to the deformation parameter characterising the deformed 
$\U_\pardeform(iso_{\k_2\dots\k_N}(N))$ CK family. 
In any case, we wish to stress here that the dimensions of
$\pardeform$ which
result from the {\it first} contraction process (that involving $\k_1$) 
depend crucially on the assumed dependence of the dimensionless $q$ (or $z$) 
on the new deformation parameter 
$\pardeform $ and of the contraction one.
For instance, if $q=\exp(\pardeform /R)$ as above, where $R$ is a radius 
(length), $[\pardeform ]=L^{1}$.
Any other dimensions assigned to $\pardeform $ necessarily include hidden 
hypotheses on the dependence of $q$ on {\it other} fundamental constants 
($\hbar$ and $c$ are necessary, for instance, in order to have 
$[\pardeform ]$=
(mass)$^{-1}$).
This is specially important because the appearance of Planck's $\hbar$, 
for instance, implies quantum considerations in the strict 
(\ie, physical) sense of the word.
These go beyond the purely mathematical deformation process, and 
should be accordingly discussed separately in any physical application of 
a $\pardeform$-deformed algebra 
(see \cite{AKRA} in connection with deformed Minkowski spaces).
 
It would be interesting to know whether the bicrossproduct structure is 
present in other cases, \ie, for $\k_1\ne 0$ different from the  $\k_1=0$ 
deformations considered here.
However, the general form of the deformation 
$\U_q(so_{\k_1\dots\k_N}(N+1))$ of the
general CK algebra in the `physical' (rather than the Cartan-Weyl) basis 
is unknown, and this precludes us for the moment from discussing this point.

\section*{Acknowledgments}

This work has been partially supported by the Spanish 
CICYT, DGICYT and DGES (grants AEN96--1669, 
PB94--1115, PR95--439 and PB95--0719).
J.A. and J.C.P.B. wish to thank the kind hospitality extended to them at DAMTP.
The support of St. John's College (J.A.) and an FPI grant from the Spanish
Ministry of Education and Science and the CSIC (J.C.P.B.) are also gratefully
acknowledged.

\renewcommand{\thesection}{Appendix:}
\setcounter{section}{0}
\section{Some bicrossproduct formulae}
\renewcommand{\theequation}{A.\arabic{equation}}

We reproduce here some of the
formulae needed in the main text, and refer to 
\cite{MajMajid} for details.
In these formulae ${\cal A}\equiv\U_\pardeform (T_N)$ and 
${\cal H}\equiv\U({so}_{\kiin} (N))$. The right action is 
$\alpha:{\cal A}\otimes{\cal H}\to {\cal A}\,,\, \alpha(a\otimes h)\equiv 
a\triangleleft h$, and the left coaction $\beta:{\cal H}\to {\cal 
A}\otimes{\cal H}\,,\, \beta(h)=h^{(1)}\otimes h^{(2)}\,,\,
(h^{(1)}\in {\cal A}\,,\,h^{(2)}\in {\cal H})$
so that $({\cal A},\alpha)$ $[({\cal H},\beta)]$ is a right ${\cal H}$-module
[left ${\cal A}$-comodule].
The compatibility conditions are
\beq
\v_{\cal A}(a\triangleleft h)=\v_{\cal A}(a)\v_{\cal H}(h)
\quad,
\label{apix}
\eeq
\beq
\Delta(a\triangleleft h)\equiv (a\triangleleft h)_{(1)}\otimes 
(a\triangleleft h)_{(2)}
=(a_{(1)}\triangleleft  h_{(1)})h_{(2)}^{\ (1)}\otimes a_{(2)}\triangleleft 
h_{(2)}^{\ (2)}\quad,
\label{apx}
\eeq
\beq
\beta(1_{\cal H})\equiv 1_{\cal H}^{(1)}\otimes 1_{\cal H}^{(2)}=
1_{\cal A}\otimes 1_{\cal H}\quad,
\label{apxi}
\eeq
\beq
\beta(hg)\equiv (hg)^{(1)}\otimes (hg)^{(2)}
=(h^{(1)}\triangleleft g_{(1)})g_{(2)}^{\ (1)}\otimes h^{(2)}g_{(2)}^{\ (2)}
\quad,\label{apxii}
\eeq
\beq
h_{(1)}^{\ (1)}(a\triangleleft  h_{(2)})\otimes h_{(1)}^{\ (2)}=
(a\triangleleft  h_{(1)})
h_{(2)}^{\ (1)}\otimes h_{(2)}^{\ (2)}
\label{apxiii}
\eeq
(subindices refer to coproduct as usual; {\it superindices} refer to the 
components of $\beta(h)$).
When they are satisfied the right-left bicrossproduct structure 
${\cal H}^\beta\bicross_\alpha {\cal A}$
on ${\cal K}\equiv{\cal H}\otimes{\cal A}$ is determined by
\beq
(h\otimes a)(g\otimes b)=hg_{(1)}\otimes (a\triangleleft  g_{(2)})b\quad,
\quad h,g\in {\cal H}
\; ; a,b\in {\cal A}\quad,
\label{apxiv}
\eeq
\beq
\Delta_{\cal K}(h\otimes a)=h_{(1)}\otimes h_{(2)}^{\ (1)}a_{(1)}
\otimes h_{(2)}^{\ (2)}\otimes a_{(2)}\quad,
\label{apxv}
\eeq
\beq
\v_{\cal K}=\v_{\cal H}\otimes\v_{\cal A}\quad,\quad
 1_{\cal K}=1_{\cal H}\otimes 1_A\quad,
\label{apxvi}
\eeq
\beq
S(h\otimes a)=(1_{\cal H}\otimes S_{\cal A}(h^{(1)}a))(S_H(h^{(2)})\otimes 
1_{\cal A})
\quad.
\label{apxvii}
\eeq

\setcounter{equation}{0}
\renewcommand{\theequation}{A'.\arabic{equation}}

Let $A,H$ be the duals of ${\cal A},{\cal H}$.
The dual `group' aspect of the above formulae imply the existence of mappings 
$\bar\alpha:A\otimes H\to H\,,\,\bar\beta:A\to A\otimes H$, dual to 
$(\beta,\alpha)$ respectively, which satisfy the conditions
\beq
\v_{\AA}(\hh\RRR \aa)=\v_{\HH}(\hh)\v_{\AA}(\aa)\quad,
\label{abpix}
\eeq
\beq
\Delta(\hh\RRR \aa)\equiv (\hh\RRR \aa)_{(1)}\otimes (\hh\RRR \aa)_{(2)}
=(\hh_{(1)}^{\ (1)}\RRR \aa_{(1)})\otimes \hh_{(1)}^{\ (2)}(\hh_{(2)}
\RRR \aa_{(2)})
\quad,
\label{abpx}
\eeq
\beq
\bbeta(1_{\HH})\equiv 1_{\HH}^{(1)}\otimes 1_{\HH}^{(2)}=
1_{\HH}\otimes 1_{\AA}\quad,
\label{abpxi}
\eeq
\beq
\bbeta(\hh\gg)\equiv (\hh\gg)^{(1)}\otimes (\hh\gg)^{(2)}
=\hh_{(1)}^{\ (1)}\gg^{(1)}\otimes \hh_{(1)}^{\ (2)}(\hh_{(2)}\RRR \gg^{(2)})
\quad,
\label{abpxii}
\eeq
\beq
\hh_{(2)}^{\ (1)}\otimes(\hh_{(1)}\RRR \aa)\hh_{(2)}^{\ (2)}=
\hh_{(1)}^{\ (1)}\otimes \hh_{(1)}^{\ (2)}(\hh_{(2)}\RRR \aa)\quad.
\label{abpxiii}
\eeq
Then, there is a left-right bicrossproduct structure 
$H_{\bar\alpha}\LR^{\bar\beta} A$ on $K=H\otimes A$ defined by
\beq
(\aa\otimes \hh)(\bb\otimes \gg)=\aa(\hh_{(1)}\RRR \bb)\otimes 
\hh_{(2)}\gg\quad,
\quad\aa,\bb\in {\AA}\;;\;\hh,\gg\in {\HH}\quad,
\label{abpxiv}
\eeq
\beq
\Delta_{\KK}(\aa\otimes \hh)=\aa_{(1)}\otimes \hh_{(1)}^{\ (1)}
\otimes \aa_{(2)}
\hh_{(1)}^{\ (2)}\otimes \hh_{(2)}\quad,\label{abpxv}
\eeq
\beq
\v_{\KK}=\v_{\AA}\otimes\v_{\HH}\quad,\quad
 1_{\KK}=1_{\AA}\otimes 1_{\HH}\quad,\label{abpxvi}
\eeq
\beq
S(\aa\otimes \hh)=(1_{\AA}\otimes S_{\HH}(\hh^{(1)}))(S_{\AA}(\aa\hh^{(2)})
\otimes 1_{\HH})
\quad.
\label{abpxvii}
\eeq


\end{document}